\documentstyle[epsf]{article}
\newtheorem{theorem}{Theorem}[section]
\newtheorem{lemma}[theorem]{Lemma}
\newtheorem{Conjecture}[theorem]{Conjecture}

\newcommand{\be}{\begin{equation}}
\newcommand{\ee}{\end{equation}}
\newcommand{\Rr}{\mbox{\bf R}}
  \title{On a conjecture of Bennewitz, and the behaviour of the
 Titchmarsh-Weyl matrix near a pole. }
\author{ B.M.Brown 
Department of Computer Science, \\University of Wales, Cardiff, 
\\PO Box 916, Cardiff CF2 3XF  U.K.
\and M. Marletta  Department of Mathematics and Computer Science, \\ University 
of Leicester, \\
 University Road, Leicester LE1 7RH  U.K.}
\begin{document}
\maketitle
\begin{abstract} For any real limit-$n$ $2n$th-order selfadjoint linear 
 differential expression on $[0,\infty)$, Titchmarsh-Weyl matrices
 $M(\lambda)$ can be defined. Two matrices of particular interest are
 the matrices $M_D(\lambda)$ and $M_N(\lambda)$ associated respectively 
 with Dirichlet and Neumann boundary conditions at $x=0$. These satisfy
 $M_D(\lambda) = -M_{N}(\lambda)^{-1}$. It is known that when these
 matrices have poles (which can only lie on the real axis) the existence 
 of valid HELP inequalities depends on their behaviour in the neighbourhood 
 of these poles. We prove a conjecture of Bennewitz and use it, together with 
 a new algorithm for computing the Laurent expansion of a Titchmarsh-Weyl 
 matrix in the neighbourhood of a pole, to investigate the existence of HELP
 inequalities for a number of differential equations which have so
 far proved awkward to analyse.
 \end{abstract}
\newcommand{\beq}{\begin{equation}}
\newcommand{\enq}{\end{equation}}
\newcommand{\s}{ {\cal{x}}}

\section{Introduction}

In a recent paper \cite{kn:lotsofus} a numerical algorithm was reported
for the computation of the Titchmarsh-Weyl $M$ matrices associated 
with the fourth order differential equation
\beq
{\cal M}[y] = ((py'')'-(sy'))' +qy=\lambda y. \label{eq:1.1}
\enq
The motivation for that work was an investigation of the HELP type integral 
inequality
\be
\left( \int_0^\infty ( p \mid y'' \mid ^2 + s \mid y ' \mid ^2 +q \mid y \mid^2) 
dx\right)^2 
\leq K \int_0^\infty \mid y \mid ^2 dx \int_0^\infty \mid {\cal M}[y] \mid ^2 
dx. 
\label{eq:1.2}
\ee
In this paper we turn our attention to higher order operators and inequalities.
We consider operators of the form
\be {\cal M}[y] := 
\sum_{j=0}^{n}(-1)^j\frac{d^j}{dx^j}\left(p_j(x)\frac{d^jy}{dx^j}\right); 
\label{eq:1.1star} \ee
as we do not wish to be concerned with quasidifferential expressions, we assume 
that
the coefficients $p_j$ are all smooth. We assume that $p_n > 0$ on $(0,\infty)$,
and we assume that ${\cal M}$ is regular at $x=0$. ${\cal M}$ will possess 
selfadjoint
realisations in $L^2[0,\infty)$: we assume that $x=\infty$ is of limit-point 
(minimal
deficiency index) type, so that these realisations depend only upon choice of 
boundary 
conditions at $x=0$. The corresponding HELP inequality
is
\be \left(\int_0^{\infty} \left\{ \sum_{j=0}^{n}p_j|y^{(j)}|^2 \right\} dx 
\right)^2
 \leq K \left(\int_0^\infty |y|^2 dx\right) \left(\int_0^\infty |{\cal 
M}[y]|^2dx\right).
\label{eq:1.2star} 
\ee
\par
Before we explain the connection between Titchmarsh-Weyl $M$-matrices and 
inequalities 
(\ref{eq:1.2}) and (\ref{eq:1.2star}), we give a very brief overview of 
$M$-matrices, starting
with the scalar case ($1\times 1$ matrices). Consider a second order 
Sturm-Liouville 
equation, say 
\be -(py')'+qy = \lambda w y, \label{eq:mmrev1} \ee 
on an interval $[0,\infty)$, with $x=0$ a regular point and $x=\infty$ a 
singular point of limit-point type.  Suppose that $y_D$ denotes the solution of 
(\ref{eq:mmrev1}) 
subject to the Dirichlet conditions
\[ y_D(0) = 0, \;\;\; py_D'(0) = 1, \]
while $y_N$ denotes the solution subject to the Neumann conditions
\[ y_N(0) = -1, \;\;\; py_N'(0) = 0. \]
Then the Dirichlet $m$-function $m_D(\lambda)$ is, for $\Im(\lambda)\neq 0$,
the unique function such that
\[ y_N(\cdot) + m_D(\lambda) y_D(\cdot) \]
is a solution of (\ref{eq:mmrev1}) square integrable over $[0,\infty)$ 
with respect to $w$. The Neumann $m$-function $m_N(\lambda)$ is defined by the 
property that
\[ y_D(\cdot) - m_N(\lambda) y_D(\cdot) \]
is a square integrable over $[0,\infty)$. The functions $m_D$ and $m_N$ are 
analytic functions 
of $\lambda$ on both the upper and lower half planes; moreover 
$m_D(\lambda)m_N(\lambda) = -1$ 
wherever both $m_D$ and $m_N$ are defined. 

These ideas were generalized by Hinton and Shaw \cite{kn:Hinton} to certain 
Hamiltonian 
systems (which include (reformulations of) higher order Sturm-Liouville systems
such as (\ref{eq:1.1star}) -- see Section \ref{section:3.2} below for details
of the fourth order case). Consider a system
\be  \left( \begin{array}{c} -v' \\ u' \end{array}\right) = S(x;\lambda)
   \left( \begin{array}{c} u \\ v \end{array}\right) \label{eq:mmham} \ee
in which $S$ is a $2n\times 2n$ symmetric matrix given in terms of real
matrices $A$ and $B$ by
\[ S(x;\lambda) = A(x) + \lambda B(x), \]
where $B$ is positive semi-definite and, in a certain sense, positive definite 
on solutions of the differential equation. The dependent variables  $v$ and $u$ 
are $n$-vector functions of $x$ and $\lambda$. Let 
\[ \left( \begin{array}{c} U_D \\ V_D \end{array}\right) \]
be the $2n\times n$ matrix solution of this ODE, partitioned into $n\times n$
blocks, subject to the initial conditions 
\[ U_D(0) = {\bf 0}, \;\;\; V_D(0) = I, \]
and let 
\[ \left( \begin{array}{c} U_N \\ V_N \end{array}\right) \]
be the $2n\times n$ matrix solution subject to initial conditions 
\[ U_N(0) = -I, \;\;\; V_N(0) = {\bf 0}. \]
 Then the Dirichlet and Neumann $M$-matrices are defined respectively by the 
requirements that
\be U_1(\cdot) := U_N(\cdot) + M_D(\lambda) U_D(\cdot), \;\;\; 
   U_2(\cdot) := U_D(\cdot) - M_N(\lambda) U_N(\cdot) \label{eq:mmrev17} \ee
be square integrable with respect to $B$ over $[0,\infty)$, for $\Im(\lambda) 
\neq 0$:
\[ \int_{0}^{\infty} U_1^*(x)B(x)U_1(x)dx < +\infty,
   \int_{0}^{\infty} U_2^*(x)B(x)U_2(x)dx < +\infty. \] 
{\bf The limit-point hypothesis at infinity ensures that $M_N$ and $M_D$ are 
uniquely 
determined by these conditions.}

We shall also make extensive use of the identity
\[ M_D(\lambda)M_N(\lambda) = -I \]
which holds wherever $M_D$ and $M_N$ are defined. For further details see 
Hinton and Shaw \cite{kn:Hinton}.

Returning to the HELP type inequality (\ref{eq:1.2}), we remark that this has 
been 
investigated by Russell \cite{kn:Russell}, while the more general form 
(\ref{eq:1.2star}) 
has been studied by Dias \cite{kn:Dias}. The form (\ref{eq:1.2}) was also 
investigated
earlier in a somewhat more restricted form by Bradley and Everitt 
\cite{kn:BradleyE} 
and Brodlie and Everitt \cite{kn:BrodlieE}. In all these investigations the 
existence of a 
valid inequality, that is a finite number $K$ in (\ref{eq:1.2}) or 
(\ref{eq:1.2star}), 
was shown to depend upon the behaviour of the Titchmarsh-Weyl matrix $M_N$ 
associated
with (\ref{eq:1.1}) or (\ref{eq:1.1star}).
 The existence of the inequality and value of the best constant is determined by 
the behaviour of the function
\be
 \Im ( \lambda^2 M_N(\lambda))  \;\;\; \mbox{(the imaginary part of $\lambda^2
 M_N(\lambda)$)} \label{eq:1.3}
\ee
for strictly complex values of the spectral parameter $\lambda$ that lie in 
the first and third quadrants of the complex plane. Indeed the existence, 
but not necessarily the value, of the best constant is determined by 
(\ref{eq:1.3}) for values such that $\mid \lambda \mid  \rightarrow 0 $.
As it is  difficult to find examples of $M$ matrices for  (\ref{eq:1.1}) or 
(\ref{eq:1.1star}) 
which are known in closed form, it is of some importance to be able to 
investigate this 
problem numerically.
\par
It is further known that when $0$ lies both in the resolvent set of the 
realisation of
${\cal M}$ generated by Neumann boundary conditions ($v(0)=0$ in the Hamiltonian 
formulation,
or $-(py'')'(0)+sy'(0)=0$, $py''(0)=0$ in the fourth order case) and also in the 
resolvent 
set of the realisation generated by Dirichlet boundary conditions ($u(0) = 0$, 
or
$y(0) = 0$, $y'(0)=0$ in the fourth order case) then the 
inequality (\ref{eq:1.2star}) fails; a necessary (though not generally 
sufficient)
condition for an inequality is that $0$ be a point of the spectrum of at least 
one of 
these two operators. However, when the Titchmarsh-Weyl matrices are meromorphic,
a little more can be said on the validity of the inequality. 
In this case Dias \cite{kn:Dias} has shown that the poles of $M_N$ occur at the 
eigenvalues of the realisation of ${\cal M}$ subject to Neumann conditions 
$v(0)=0$ and 
the poles of $M_D$ occur at the eigenvalues of ${\cal M}$ subject to Dirichlet 
conditions 
$u(0)=0$.

If  $\mu$ is a pole of an $M$-matrix then $\mu$ is simple \cite{kn:Hinton} and 
the 
$M$-matrix has an expansion
\be
M(\lambda) = \frac{\sigma_{-1}}{ \lambda -\mu} + \sigma_0 +...
\ee
where $\sigma_{-1}$ is real and is called the residue matrix of $M$ at $\mu$.
We shall need the concepts of  {\bf Neumann } and {\bf Dirichlet }
translates. In (\ref{eq:1.1star}) we say that $\mu$ is a {\bf Neumann } 
translate 
if it is an eigenvalue of 
\be
{\cal M}[y] = \mu y \label{eq:1.4}
\ee
with Neumann condition given in the Hamiltonian form by $v(0) = 0$.
In a similar way we say that $\mu$ is a {\bf Dirichlet } translate if  it is 
an eigenvalue of (\ref{eq:1.4}),  but this time with initial conditions 
$u(0)=0$. The residue matrix of $M_N$ associated with a Neumann 
translate we denote by $\sigma_N$, and the residue matrix of $M_D$ subject
to a Dirichlet translate we denote by $\sigma_D$.

Suppose that  $\mu_N, \mu_D$ are Neumann and Dirichlet translates 
respectively. It  is shown in \cite{kn:Dias} that a valid inequality  
will be found if the differential expression (\ref{eq:1.1star}) is replaced 
by either
\[ 
{\cal M}_N = {\cal M}-\mu_N \]
or 
\[
{\cal M}_D = {\cal M}-\mu_D  \]
provided either of the associated residue matrices $\sigma_N, \sigma_D$ 
is of full rank. The result for higher order HELP type inequalities is 
somewhat weaker than that for the second order problem since it is shown in
\cite{kn:Everitt}
for the second order classical HELP inequality that when the $m$ function 
is meromorphic then  a valid  inequality  exists if and only if  $0$ is 
an eigenvalue of either the Neumann or the Dirichlet problem associated 
with that expression.   In an attempt to strengthen this result for higher 
order operators Bennewitz (private communication, 1995) has proposed the 
conjecture that  provided 
\be
{\rm rank } (\sigma_N) + {\rm rank } (\sigma_D) =n  \label{eq:1.6}
\ee
(half  the order of the differential expression) then a valid inequality 
will be found. In this paper we shall prove this conjecture for the general
even order HELP inequality.
\par
As we remarked above  the existence of a valid inequality is determined by 
(\ref{eq:1.3}) as $\mid \lambda \mid \rightarrow 0$ and when $M_N$ is 
meromorphic it must also have a pole at $0$. Thus in order to 
investigate numerically the existence of a valid inequality  we must compute
the associated  residue matrices $\sigma_N, \sigma_D$. In \cite{kn:lotsofus} 
it was noted that this was a difficult numerical problem and in 
section \ref{section:algorithms} we report on  
some new algorithms to solve it.

In section \ref{section:numerics} we apply our work to some equations
to determine whether or not they are likely to possess associated
HELP inequalities.

\section{The Bennewitz Conjecture}
In order to simplify the algebra we shall assume that at least
one of the matrices $M_{D}(\lambda)$, $M_N(\lambda)$ has a pole
at the origin $\lambda=0$. The Bennewitz Conjecture is as follows.
\begin{Conjecture} (Bennewitz) Suppose that $M_N$ (and hence $M_D$)
is meromorphic and that
\[ \mbox{rank}(\mbox{Res}(M_D,0)) + 
   \mbox{rank}(\mbox{Res}(M_N,0)) = n. \]
Then there is a valid HELP inequality associated with the differential
expression.
\end{Conjecture}
\subsection{Proof of the Bennewitz Conjecture}
In order to prove this result we shall require the following
lemmas.
\begin{lemma} The Titchmarsh-Weyl matrices $M_D$ and $M_N$ are
symmetric matrices and are also Nevanlinna functions, in the sense 
that the matrices $\Im(M_D(\lambda))$ and $\Im(M_N(\lambda))$ are 
positive definite for $\Im(\lambda)>0$. 
\end{lemma}
For a proof of this result see Hinton and Shaw \cite{kn:Hinton}.
\begin{lemma} \label{lemma:2}
A necessary and sufficient condition for the existence
of a HELP inequality is that there exist numbers $\theta_{+} \in (0,\pi/2)$,
$\theta_{-}\in (\pi,3\pi/2)$,  $\rho_{+}>0$ and $\rho_{-}>0$ such that 
\[ \Im(-\lambda^2 M_N(\lambda)) > 0 \;\;\; \forall |\lambda|\in (0,\rho_{+}),
 \; \mbox{arg}(\lambda )\in [\theta_{+},\pi/2) \]
and
\[ \Im(\lambda^2 M_N(\lambda)) > 0 \;\;\; \forall |\lambda|\in (0,\rho_{-}),
 \; \mbox{arg}(\lambda )\in [\theta_{-},3\pi/2). \]
\end{lemma}
This result is proved in Dias' thesis \cite{kn:Dias}, and in a different
form for the case $n=2$ in Russell \cite{kn:Russell}.
\vspace{2mm}

\noindent {\bf Notes}
\par
\begin{enumerate}
\item In Lemma \ref{lemma:2} and throughout the rest of this paper, we 
 follow the usual convention for matrices that relations of the form
 `$>0$' and `$<0$' indicate positive definiteness and negative definiteness
 respectively.
\item We could state this lemma in an equivalent form in which there would
be just one number $\rho>0$ equivalent to $\min(\rho_{+},\rho_{-})$. 
However for the proof that follows this form is marginally more
convenient.
\end{enumerate}
\begin{lemma} \label{lemma:mm3} 
Suppose that $M_D$ or $M_N$ has a simple pole at $\lambda = 0$ with a 
Laurent expansion
\[ \frac{1}{\lambda}M_{-1} + M_0 + \lambda M_1 + \lambda^2 M_2 + \cdots.\]
Then all the coefficients occurring in this expansion are real symmetric
matrices.
\end{lemma}
\noindent {\bf Proof}\, The symmetry of the coefficients follows from the
symmetry of $M_N$ and $M_D$, see Hinton and Shaw \cite{kn:Hinton}. For the
rest of the proof we concentrate on $M_N$: the proof for $M_D$ is
similar.

>From (\ref{eq:mmrev17}) with $x=0$ it is clear that 
\[ M_N(\lambda) = U_{2}(0;\lambda), \]
where $\left(\begin{array}{c} U_2 \\ V_2 \end{array}\right)$ is a
`square integrable' solution of the Hamiltonian system (more
precisely, is a solution for which $U_2$ is square integrable in the
sense of Hinton and Shaw). Also, it can be shown from (\ref{eq:mmrev17}) that
\[ V_2(x;\lambda) = V_D(x;\lambda) - M_N(\lambda)V_N(x;\lambda), \]
whence setting $x=0$ gives 
\[ V_2(0;\lambda) = I. \]
Thus
\be M_N(\lambda) = U_{2}V_{2}^{-1}(0;\lambda). \label{eq:mmrev18} \ee
Now it is easy to see that if 
$\left(\begin{array}{c} U_2(x;\mu) \\ V_2(x;\mu) \end{array}\right)$
is a square integrable solution for $\lambda=\mu$ then
$\left(\begin{array}{c} \overline{U_2(x;\mu)} \\ \overline{V_2(x;\mu)} 
\end{array}\right)$
is a square integrable solution for $\lambda = \overline{\mu}$. Since
we are concerned with problems of limit-point type, the square integrable 
solution for any $\Im(\lambda) \neq 0$ is unique up to postmultiplication by an
invertible constant matrix. Any such matrix cancels out upon taking the
combination $U_2V_2^{-1}$, and hence
\[ \overline{U_{2}V_{2}^{-1}(0;\mu)} = U_{2}V_{2}^{-1}(0;\overline{\mu}). \]
Thus from (\ref{eq:mmrev18}), 
\[ M_N(\overline{\lambda}) = \overline{M_N(\lambda)}. \]
This implies that the coefficients in the Laurent expansion of $M_N$ about
$\lambda = 0$ are real matrices. \hfill $\Box$
\vspace{2mm}

\noindent {\bf Proof of the Bennewitz Conjecture} 
\noindent Under the hypotheses of the conjecture, we can expand
$M_D$ and $M_N$ in Laurent series about the point $\lambda=0$,
\begin{eqnarray}
M_{D}(\lambda) & = & \lambda^{-1}M_{-1} + M_{0} + \lambda M_{1} + O(\lambda^2),
 \label{eq:mmMDrep} \\
M_{N}(\lambda) & = & \lambda^{-1}\hat{M}_{-1} + \hat{M}_{0} + 
 \lambda \hat{M}_{1} + O(\lambda^2), \label{eq:mmMNrep}
\end{eqnarray}
the expansions being valid in a neighbourhood of $\lambda=0$. By Lemma 
\ref{lemma:mm3}
the coefficients in these expansions are real symmetric matrices. We see from 
the hypothesis
of the conjecture that 
\[ \mbox{rank}(M_{-1}) + \mbox{rank}(\hat{M}_{-1}) = n, \]
so let $\mbox{rank}(M_{-1}) = r$, $\mbox{rank}(\hat{M}_{-1})=n-r$.
We also know that $M_{D}M_{N} = M_{N}M_{D} = -I$ for all
$\Im(\lambda)>0$ (see Hinton and Shaw \cite{kn:Hinton}). Multiplying the 
Laurent expansions, we obtain the following conditions.
\be
M_{-1}\hat{M}_{-1} = \hat{M}_{-1}M_{-1} = 0, \label{eq:mmnull}
\ee
\be
M_{-1}\hat{M}_{0} + M_{0}\hat{M}_{-1} = \hat{M}_{-1}M_{0} + 
 \hat{M}_{0}M_{-1}= 0, \label{eq:mmnull3}
\ee
\be 
M_{0}\hat{M}_0 + M_{-1}\hat{M}_{1} + M_{1}\hat{M}_{-1} = 
\hat{M}_0M_0 + \hat{M}_{-1}M_1 + \hat{M}_1M_{-1} = -I.
\ee 
Thus the columns of $M_{-1}$ are orthogonal to the columns of
$\hat{M}_{-1}$. We can therefore choose an orthonormal basis
of $\Rr^n$, say $\{ v_{1},\ldots, v_n\}$, such that 
\[ \mbox{Span}\{v_1,\ldots,v_r\} = \mbox{Span of the columns of $M_{-1}$}, \]
\[ \mbox{Span}\{v_{r+1},\ldots,v_n\} = 
    \mbox{Span of the columns of $\hat{M}_{-1}$}. \]
Let $V$ be the $n\times r$ matrix with columns $v_1,\ldots, v_r$,
and let $\hat{V}$ be the $n\times (n-r)$ matrix with columns
$v_{r+1},\ldots,v_n$. We make the following observations. 
\begin{itemize} 
\item From (\ref{eq:mmnull}),
\be M_{-1}\hat{V} = 0; \;\; \hat{V}^TM_{-1} = 0; \;\;
    \hat{M}_{-1}V = 0; \;\; V^T\hat{M}_{-1} = 0.
\label{eq:mmnull2} 
\ee
\item Any vector $v\in \Rr^n$ can be written in the form
\be v = V\alpha + \hat{V}\hat{\alpha}, \label{eq:mmvrep} \ee
where $\alpha\in \Rr^{r}$ and $\hat{\alpha}\in \Rr^{n-r}$.
\end{itemize}
Using (\ref{eq:mmvrep}) 
%together with the fact that $V$ and $\hat{V}$ are real matrices, 
we have
\begin{eqnarray}
v^T M_N v & = & \alpha^T V^T M_N V \alpha + 
            \hat{\alpha}^T \hat{V}^TM_N\hat{V}\hat{\alpha} \nonumber \\
& + & 2\alpha^T V^TM_N \hat{V}\hat{\alpha}, \label{eq:mmvMv}
\end{eqnarray}
%\begin{eqnarray}
%v^T \Im (M_N) v & = & \alpha^T \Im (V^T M_N V) \alpha + 
%            \hat{\alpha}^T \Im (\hat{V}^TM_N\hat{V})\hat{\alpha} \nonumber \\
%& + & 2\alpha^T \Im (V^TM_N \hat{V})\hat{\alpha}, \label{eq:mmvMv}
%\end{eqnarray}
where we have used the symmetry of $M_N$ to simplify the last term.
We now use the Laurent expansion (\ref{eq:mmMNrep}) together with
the conditions (\ref{eq:mmnull2}) to simplify this expression. 
We observe that 
\begin{eqnarray}
  V^T M_N V & = & V^T\left(\lambda^{-1}\hat{M}_{-1} + \hat{M}_0 + 
  \lambda \hat{M}_1 + O(\lambda^2) \right) V \nonumber \\
 & = & V^T \hat{M}_0 V + \lambda V^T \hat{M}_1 V + O(\lambda^2).
 \nonumber
\end{eqnarray}
Now combining (\ref{eq:mmnull}) and (\ref{eq:mmnull3}) we obtain
$M_{-1}\hat{M}_0M_{-1} = 0$, which implies that $V^T \hat{M}_0 V = 0$.
Thus
\be
  V^T M_N V  =  \lambda V^T \hat{M}_1 V + O(\lambda^2). \label{eq:mmterm1}
\ee
At this stage it is not clear that the matrix $V^T \hat{M}_1 V$ is
of full rank, so the $O(\lambda^2)$ terms might not be negligible.
We shall show shortly that in fact $V^T \hat{M}_1 V$ is of full rank.
The next term we need to examine is $\hat{V}^T M_N \hat{V}$. This
is much easier to deal with; the Laurent expansion (\ref{eq:mmMNrep})
immediately gives
\be 
 \hat{V}^T M_N \hat{V} = \frac{1}{\lambda} \hat{V}^T \hat{M}_{-1} \hat{V}
 + O(1), \label{eq:mmterm2}
\ee
and the leading order term here is of full rank since the span
of the columns of $\hat{V}$ is the same as the span of the columns
of $\hat{M}_{-1}$, and $\hat{M}_{-1}$ is symmetric. Finally, we treat the
term $V^T M_N \hat{V}$. Using (\ref{eq:mmMNrep}) and (\ref{eq:mmnull2}) we 
obtain
\be V^T M_N \hat{V} = V^T \hat{M}_0 \hat{V} + \lambda V^T \hat{M}_1 \hat{V}
 + O(\lambda^2). \label{eq:mmterm3}
\ee
Once more, since we know little about the ranks of the coefficients
in this expansion, it is not clear that the $O(\lambda^2)$ terms are
negligible, a point which will have to be borne in mind later on.
Substituting (\ref{eq:mmterm1}), (\ref{eq:mmterm2}) and (\ref{eq:mmterm3})
back into (\ref{eq:mmvMv}) we obtain 
\begin{eqnarray} v^T M_N v & = & \alpha^T 
 \left(\lambda V^T \hat{M}_1 V + O(\lambda^2) \right) \alpha \nonumber \\
 & + & \hat{\alpha}^T  \left(\frac{1}{\lambda} \hat{V}^T \hat{M}_{-1} \hat{V} 
+ O(1) \right)\hat{\alpha} \nonumber \\
 & + & 2\alpha^T\left(V^T \hat{M}_0 \hat{V} + \lambda V^T \hat{M}_1 \hat{V}
 + O(\lambda^2)\right) \hat{\alpha}. \label{eq:mmeq5} 
\end{eqnarray} 
>From the Nevanlinna property of $M_N$ we know that we must have
$v^T \Im(M_N) v > 0$ for all $\Im(\lambda) > 0$, for all non-zero 
$v\in \Rr^n$. Choosing $\hat{\alpha}=0$ we see that this implies, in 
particular, that
\be \alpha^T \Im \left(\lambda V^T \hat{M}_1 V + O(\lambda^2) \right) \alpha
 > 0 \label{eq:mmineq4} \ee
for all non-zero $\alpha \in \Rr^r$. Now suppose that $V^T \hat{M}_1 V$ is
not of full rank. Then we can choose a non-zero $\alpha$ such that
\[ \alpha^T (V^T\hat{M}_1 V)\alpha = 0. \]
Suppose now that with this choice of $\alpha$ we have an expansion of
the form
\[ v^T \Im(M_N) v = \alpha^T(V^T \hat{M}_{p}V) \alpha \Im(\lambda^p) + O(
 \Im(\lambda^{p+1})). \]
where $p>1$ and the coefficient $\alpha^T(V^T \hat{M}_{p}V) \alpha $ is
non-zero. Because $p>1$, we know that $\Im(\lambda^p)$ is not of one
sign on the upper half plane. Thus $ v^T \Im(M_N) v $ cannot be of one
sign in the upper half plane, contradicting the Nevanlinna property of
$M_N$. We have thus established that the matrix $V^T \hat{M}_1 V$ is
of full rank; also, from (\ref{eq:mmineq4}), we have therefore established
that it is positive definite:
\be V^T \hat{M}_1 V > 0. \label{eq:mmposdef1} \ee
Next we choose $\alpha = 0$, $\hat{\alpha} \neq 0$ in (\ref{eq:mmeq5}):
since $\Im(\frac{1}{\lambda})<0$ when $\Im(\lambda)>0$,  we see that 
$\hat{V}^T \hat{M}_{-1} \hat{V}$ is negative definite, i.e.
\be \hat{V}^T \hat{M}_{-1} \hat{V} < 0. \label{eq:mmposdef2} \ee
The two results (\ref{eq:mmposdef1}) and (\ref{eq:mmposdef2}) -- together
with (\ref{eq:mmeq5}) -- imply the Bennewitz conjecture, as we shall show
in the remainder of the proof. 

>From Lemma \ref{lemma:2}, we first need to show that there exists 
$\theta_+ \in (0,\pi/2)$ and $\rho_{+}>0$ such that for $|\lambda| \in 
(0,\rho_{+})$ 
and $\mbox{arg}(\lambda)\in [\theta_+,\pi/2)$, the matrix 
$\Im(-\lambda^2 M_N(\lambda))$ is positive definite.

With $v = V\alpha + \hat{V}\hat{\alpha}$ as before, (\ref{eq:mmeq5}) gives
\begin{eqnarray} v^T \Im(-\lambda^2 M_N) v & = & \alpha^T 
\Im \left(-\lambda^3 V^T \hat{M}_1 V + O(\lambda^4) \right) \alpha \nonumber \\
 & + & \hat{\alpha}^T \Im \left(-\lambda \hat{V}^T \hat{M}_{-1} \hat{V} 
+ O(\lambda^2) \right)\hat{\alpha} \nonumber \\
 & + & 2\alpha^T\Im\left(-\lambda^2 V^T \hat{M}_0 \hat{V} - \lambda^3 V^T 
\hat{M}_1 \hat{V}
 + O(\lambda^4)\right) \hat{\alpha}. \label{eq:mmeq6} 
\end{eqnarray} 

Let $\lambda = \rho \mbox{e}^{i\theta}$. Fix a number $\theta_1 \in 
(\pi/3,\pi/2)$.
Then $-\sin(3\theta) > 0$ for $\theta \in [\theta_1,\pi/2]$.
Thus (\ref{eq:mmposdef1}) implies
that there is a constant $\omega_1 > 0$ and a number $\rho_{1} > 0 $ 
such that for all $\rho \in (0,\rho_1)$ and $\theta\in [\theta_1,\pi/2]$ the 
first 
term in (\ref{eq:mmeq6}) satisfies
\be \alpha^T \Im \left(-\lambda^3 V^T \hat{M}_1 V + O(\lambda^4) \right) \alpha
 > \omega_1\rho^3 \| \alpha \|^2. \label{eq:mmineq5}
\ee
Similarly, fix $\theta_2 \in (0,\pi/2)$. Then $\sin\theta >0$ for $\theta 
\in [\theta_2,\pi-\theta_2]$. Thus
(\ref{eq:mmposdef2}) implies that there is a constant $\omega_2>0$ and 
a number $\rho_2 > 0 $ such that for all $\rho \in (0,\rho_2)$ and 
$\theta \in [\theta_2,\pi-\theta_2]$ the second 
term in (\ref{eq:mmeq6}) satisfies
\be
\hat{\alpha}^T \Im \left(-\lambda \hat{V}^T \hat{M}_{-1} \hat{V} 
+ O(\lambda^2) \right)\hat{\alpha} > \omega_2 \rho \| \hat{\alpha} \|^2.
\label{eq:mmineq6}
\ee
We now deal with the last term in (\ref{eq:mmeq6}).
Clearly there exist positive constants $C_1$ and $C_2$ and $r$ such that 
for all $\rho=|\lambda| \in (0,r)$ and $\theta = \mbox{arg}(\lambda) \in 
(0,\pi)$,
\[ \left| 2\alpha^T\Im\left(-\lambda^2 V^T \hat{M}_0 \hat{V} - 
                     \lambda^3 V^T \hat{M}_1 \hat{V} + 
                     O(\lambda^4)\right) \hat{\alpha} \right|
   \leq \left(C_1 \rho^2 |\sin(2\theta)| + C_2 \rho^3\right) 
        \| \alpha \| \|\hat{\alpha}\|.
\]
We bound the second part using Young's inequality:
\[ C_2 \rho^3 \| \alpha \| \|\hat{\alpha}\|
   \leq \frac{1}{2}C_2\rho
   \left\{ \rho^3 \| \alpha \|^2 + \rho \|\hat{\alpha}\|^2\right\}.
\]
We also bound the first part using Young's inequality:
\[ C_1 \rho^2 |\sin(2\theta)|\| \alpha \| \|\hat{\alpha}\|
   \leq \frac{1}{2}C_1 |\sin(2\theta)| \left\{ \rho^3 \| \alpha\|^2 +
   \rho \| \hat{\alpha} \|^2 \right\}.
\]
Combining these inequalities we obtain
\[ \left|2\alpha^T\Im\left(-\lambda^2 V^T \hat{M}_0 \hat{V} - 
                     \lambda^3 V^T \hat{M}_1 \hat{V} + 
                     O(\lambda^4)\right) \hat{\alpha}\right|
 \hspace{2.9in} \]
\be \hspace{1.5in}\leq \| \alpha\|^2 \left\{ \frac{C_1}{2}\rho^3 |\sin(2\theta)| 
+
      \frac{C_2}{2}\rho^4\right\} 
 + \| \hat{\alpha}\|^2 \left\{ \frac{C_1}{2}\rho |\sin(2\theta)| + 
\frac{C_2}{2}\rho^2\right\}. 
\label{eq:mmineq8} 
\ee
We now combine (\ref{eq:mmineq5}), (\ref{eq:mmineq6}) and (\ref{eq:mmineq8}).
Let $\theta^{*} = \max(\theta_1,\theta_2)\in (0,\pi/2)$ and let
$\rho^{*} = \min(\rho_1,\rho_2,r)$. Then for all $\rho = 
|\lambda| \in (0,\rho^{*})$ and for all $\theta = \mbox{arg}(\lambda) 
\in [\theta^{*},\pi/2]$ we have, from (\ref{eq:mmeq6}),
\be v^T \Im(-\lambda^2 M_N)v > 
\left(\omega_1 - \frac{C_1}{2}|\sin(2\theta)|-\frac{C_2}{2}\rho\right)\rho^3
 \| \alpha\|^2 + \left(\omega_2 - 
\frac{C_1}{2}|\sin(2\theta)|-\frac{C_2}{2}\rho\right)
 \rho \| \hat{\alpha} \|^2 
\ee
Choosing $\theta_{+} \in [\theta^{*},\pi/2)$ sufficiently close to
$\pi/2$ (to make $|\sin(2\theta)|$ small) and choosing $\rho_{+}
\in (0,\rho^{*}]$ sufficiently small, we can ensure that
\be v^T \Im(-\lambda^2 M_N)v > 
 \frac{\omega_1}{2}\rho^3 \| \alpha\|^2 + \frac{\omega_2}{2} \rho \| 
 \hat{\alpha} \|^2 \label{eq:mmfix1}
\ee
for all $\rho = |\lambda| \in (0,\rho_+)$ and $\theta =\mbox{arg}(\lambda)
 \in [\theta_{+},\pi/2]$. This implies that $\Im(-\lambda^2M_N(\lambda))
$ is positive definite for all such $\lambda$, which deals with the first
condition in Lemma \ref{lemma:2}.

To verify the second condition in Lemma \ref{lemma:2} we must show that there 
exists $\theta_{-} \in (\pi,3\pi/2)$ and $\rho_{-} >0$
such that $\Im(-\lambda^2 M_N(\lambda))$ is negative definite for
$|\lambda| \in (0,\rho_{-})$ and $\mbox{arg}(\lambda) \in [\theta_{-},3\pi/2)$. 
Looking back at the proof above, it is clear that there are only two changes
to the reasoning. First, we need to replace (\ref{eq:mmineq5}) by a result
of the form
\be \alpha^T \Im \left(-\lambda^3 V^T \hat{M}_1 V + O(\lambda^4) \right) \alpha
 < -\omega_3\rho^3 \| \alpha \|^2, \;\;\; (\mbox{$\omega_3 >0$ constant}),
\label{eq:mmineq5c}
\ee
to hold for all $|\lambda| \in (0,\rho_3)$ and $\mbox{arg}(\lambda) \in
 [\pi + \theta_3, 3\pi/2]$ for some $\rho_3 > 0$ and $\theta_3 \in (0,\pi/2)$.
This we can do because we can choose $\theta_3 \in (0,\pi/2)$ such that
$ -\sin(3\theta) < 0 $ for all $\theta \in [\pi + \theta_3, 3\pi/2]$.
Secondly, we need to replace (\ref{eq:mmineq6}) by a result of the form
\be
\hat{\alpha}^T \Im \left(-\lambda \hat{V}^T \hat{M}_{-1} \hat{V} 
+ O(\lambda^2) \right)\hat{\alpha} < -\omega_4 \rho \| \hat{\alpha} \|^2,
\;\;\; (\mbox{$\omega_4>0$ constant}),
\label{eq:mmineq6c}
\ee
to hold for all $|\lambda| \in (0,\rho_4)$ and $\mbox{arg}(\lambda) \in
 [\pi + \theta_4, 2\pi-\theta_4]$ for some $\rho_4 > 0$ and $\theta_4 \in 
(0,\pi/2)$.
This we can do because we can choose $\theta_4 \in (0,\pi/2)$ such that
$ \sin(\theta) < 0 $ for all $\theta \in [\pi + \theta_4, 2\pi -\theta_4]$.
Choosing $\theta_{-} \geq \pi + \max(\theta_3,\theta_4)$ sufficiently close
to $3\pi/2$ (to make $|\sin(2\theta)|$ sufficiently small, as reasoned for
(\ref{eq:mmfix1})) and choosing $\rho_{-}\in (0,\min(\rho_3,\rho_4)]$ 
sufficiently
small, one obtains an inequality of the form
\be v^T \Im(-\lambda^2 M_N)v < 
 -\frac{\omega_3}{2}\rho^3 \| \alpha\|^2 - \frac{\omega_4}{2} \rho \| 
 \hat{\alpha} \|^2 
\ee
for all $\rho = |\lambda| \in (0,\rho_{-})$ and $\theta =\mbox{arg}(\lambda)
 \in [\theta_{-},3\pi/2]$. This implies that $\Im(-\lambda^2M_N(\lambda))
$ is negative definite for all such $\lambda$. Both the conditions in Lemma
\ref{lemma:2} have now been verified, and our proof is complete. $\Box$

\subsection{A note on the converse of the Bennewitz Conjecture}
It seems appropriate to indicate here why we have been unable to prove
the converse of the Bennewitz Conjecture: namely, that if
\[ \mbox{rank}(\mbox{Res}(M_D,0)) + 
   \mbox{rank}(\mbox{Res}(M_N,0)) = n - d < n \]
then there is no HELP inequality associated with the differential
operator. The key to proving such a converse would be Lemma 
\ref{lemma:2}, which is an if-and-only-if result. Following the
notation and proof of the previous section, suppose that $V$ is an $n\times r$ 
matrix whose columns are $r$ orthonormal vectors spanning the column space of
the matrix $M_{-1}$ and let $\hat{V}$ be an $n\times (n-r)$ matrix whose 
first $n-r-d$ columns are orthonormal vectors spanning the column space of 
$\hat{M}_{-1}$ and whose last $d$ column vectors are chosen so that the 
columns of $\hat{V}$ are orthonormal. Then one can show that the columns
of $V$ and of $\hat{V}$ form an orthonormal basis of $\Rr^n$, as before;
equations (\ref{eq:mmvrep}) and (\ref{eq:mmeq5}) can be shown still to
hold. Looking for a failure in the first condition of Lemma \ref{lemma:2}, we
seek a vector $v$ such that for all $\theta_{+}\in (0,\pi/2)$ and $\rho_{+}>0$,
the inequality
\be v^T \Im(-\lambda^2 M_N(\lambda)) v > 0 \label{eq:mmeqr1} \ee
fails for some $\arg(\lambda) \in [\theta_{+},\pi/2)$ and $|\lambda| 
\in (0,\rho_{+})$. It seems reasonable to look for such a vector $v$ in
that part of $\Rr^n$ which is not spanned by the columns of $M_{-1}$ and
$\hat{M}_{-1}$: to this end we must have $\alpha = 0$ in (\ref{eq:mmvrep}),
and we must also choose $\hat{\alpha}$ such that 
$\hat{M}_{-1}\hat{V}\hat{\alpha} = 0$.
With these two conditions, (\ref{eq:mmeq5}) becomes
\be v^T \Im(M_N)v = \hat{\alpha}^T \Im(\hat{V}^T\hat{M}_0\hat{V} + 
 \lambda \hat{V}^T\hat{M}_1 \hat{V} + \ldots)\hat{\alpha}, \label{eq:mmeqr3} \ee
while similarly
\be v^T \Im(-\lambda^2 M_N)v = \hat{\alpha}^T 
\Im(-\lambda^2\hat{V}^T\hat{M}_0\hat{V} - 
 \lambda^3 \hat{V}^T\hat{M}_1 \hat{V} + \ldots)\hat{\alpha}. \label{eq:mmeqr2} 
\ee
Recall that all quantities in these equations, other than $\lambda$, are real.
Let $\lambda = r\mbox{e}^{i\theta}$. Then $\Im(-\lambda^2) = -r^2\sin(2\theta)$,
which is negative for $\theta = \pi/2 - \epsilon$, for small $\epsilon$. Thus
combining (\ref{eq:mmeqr1}) and (\ref{eq:mmeqr2}), a necessary condition for a 
HELP inequality to hold is that
\be \hat{\alpha}^T\hat{V}^T\hat{M}_0\hat{V}\alpha = 0, \label{eq:mmeqr6} \ee
whence (\ref{eq:mmeqr3}) and (\ref{eq:mmeqr2}) become
\be v^T \Im(M_N)v = \hat{\alpha}^T \Im( 
 \lambda \hat{V}^T\hat{M}_1 \hat{V} + \ldots)\hat{\alpha}, \label{eq:mmeqr4} \ee
\be v^T \Im(-\lambda^2 M_N)v = \hat{\alpha}^T \Im(- \lambda^3 \hat{V}^T\hat{M}_1 
\hat{V} 
 + \ldots)\hat{\alpha}. \label{eq:mmeqr5} \ee
The Nevanlinna condition $\Im(M_N)>0$ for $\Im(\lambda)>0$ implies that the 
leading 
term on the right hand side of (\ref{eq:mmeqr4}) is strictly positive; the 
leading
term on the right hand side of (\ref{eq:mmeqr5}) is then strictly positive for
$\arg(\lambda)\in(\pi/3,\pi/2]$, which certainly does not preclude the existence
of a HELP inequality: indeed, if it were true for all $v$ and not just those
outside the column span of $M_{-1}$ and $\hat{M}_{-1}$, it would say that a
HELP inequality definitely held. This suggests that we ought to try to prove 
that
(\ref{eq:mmeqr6}) must fail, but we have so far been unable to do this.

(Note finally that no new information is obtained by looking for a failure
in the second condition of Lemma \ref{lemma:2}). 

\section{Computing $M_N(\lambda)$ and $M_D(\lambda)$ near a pole by a 
change of variables}
\label{section:algorithms}
We now turn our attention to the problem of computing the residue matrices
of $M_N$ and $M_D$ near a pole. We shall assume once more that the pole
is at $\lambda=0$. Also, since the numerics are the same for $M_D$ as
they are for $M_N$, we shall consider a more unified problem: that of
computing the residue matrix of an arbitrary Titchmarsh-Weyl matrix $M(\lambda)$ 
having a pole at $\lambda = 0$.
 
We describe the solution of the problem in four steps. In the first of these, 
we define a new matrix $\Psi$ and show that it is
well-behaved near the pole of $M$. In the second, we explain how $\Psi$
can be computed by integrating an initial value problem. In the 
third part we explain a simple extrapolation procedure which we 
use to determine the Taylor expansion of $\Psi$ near the pole of
$M$; finally we show how the Laurent expansion of $M$ may be
recovered from the Taylor expansion of $\Psi$.

\subsection{The transformation to `safe' variables -- the matrix $\Psi$}
Suppose that the Titchmarsh-Weyl matrix $M(\lambda)$ has a pole at
$\lambda = 0$. We know that such a pole must be simple \cite{kn:Hinton}, but we 
also
know that any attempt to compute $M$ directly, by the methods we
described in our previous work \cite{kn:lotsofus}, is likely to yield inaccurate
results when $|\lambda |$ is small. In an attempt to circumvent
this difficulty we shall define a new variable $\Psi$ by
\be \Psi := (\alpha I + M^{-1})^{-1}, \label{eq:mm1} \ee
where $\alpha$ is a complex constant to be chosen for convenience.

The reason for removing the singularity in this way, rather than
by using $\lambda M(\lambda)$ as a new variable, lies in the
need to approximate whichever variable is chosen by solving an
initial value problem. For $\Psi$ the resulting Riccati-type ODE
(\ref{eq:mmrev3}) is not singular as $\lambda\rightarrow 0$.
For $\lambda M(\lambda)$, on the other hand, the corresponding
Riccati equation has a $\lambda^{-1}$ singularity occurring
in the quadratic term on the right hand side.

In order to show that $\alpha$ may be chosen so that $\Psi$ has
a removable singularity at $\lambda=0$ we need to consider two
different cases separately. The first is the case where $M^{-1}$
is well-behaved at $\lambda=0$; the second is the case where
$M^{-1}$ also has a pole at $\lambda=0$. The second of these
two cases seemed, initially, the more pathological, since it
includes the case in which the Sturm-Liouville problem with Neumann
boundary conditions shares an eigenvalue with the same problem with 
Dirichlet boundary conditions: however Lemma \ref{lemma:new} below gives
a whole class of problems for which this always happens.

\begin{lemma}
Suppose that $M^{-1}$ has a removable singularity at $\lambda = 0$. 
Then there exists a choice of $\alpha$ such that $\Psi$ has a 
removable singularity at $\lambda=0$.
\end{lemma}
\noindent {\bf Proof}\, Let $\tilde{M}$ denote 
$\lim_{\lambda\rightarrow 0}M^{-1}(\lambda)$; this exists by hypothesis.
Clearly $\Psi$ will have a removable singularity at $\lambda=0$ 
if and only if $-\alpha$ is not an eigenvalue of $\tilde{M}$. Thus
the result is proved. $\Box$

\begin{lemma}
Suppose that $M^{-1}(\lambda)$ and $M(\lambda)$ both have poles
at $\lambda=0$. Then $\Psi$ has a removable singularity at $\lambda=0$ for any
$\alpha$ with $\Im(\alpha) \neq 0$, and for all but finitely many real $\alpha$.
Also, in the case $n=2$, the matrix $\Psi$ has a removable singularity
at $\lambda=0$ for any non-zero $\alpha$.
\end{lemma}
\noindent {\bf Proof}\, 
We examine the $2\times 2$ case first. For any $2\times 2$ matrix
\[ B = \left(\begin{array}{cc} b_{1,1} & b_{1,2} \\ b_{2,1} & b_{2,2} 
 \end{array}\right) \]
let $B^A$ denote the matrix of minors of $B$, so that
\[ B^A = \left(\begin{array}{cc} b_{2,2} & -b_{1,2} \\ -b_{2,1} & b_{1,1} 
 \end{array}\right) \]
and $BB^{A} = B^{A}B = \det B$. With this notation it is clear that
\[ \Psi = \frac{\alpha I + (M^{-1})^{A}}{\det(\alpha I + M^{-1})}. \]
Expanding the determinant, we get
\be \Psi = \frac{\alpha I + (M^{-1})^{A}}{\alpha^2 + 
                \alpha \mbox{trace}(M^{-1}) + \det(M^{-1})}.
\label{eq:mm2} \ee
We also know that
\[ M^{-1} = \frac{M^A}{\det M}. \]

Since $M$ has a simple pole, so does $M^{A}$ (n.b. this step does not 
generalise to the case of $n\times n$ matrices). Also, $M^{-1}$ has a
simple pole by hypothesis. 
%Thus $\det M $ must have a removable singularity. 
Returning to (\ref{eq:mm2}), we can see that $\alpha I + (M^{-1})^A$
has a simple pole. Thus to get a removable singularity for $\Psi$, for any
non-zero $\alpha$, it suffices to show that $\mbox{trace}(M^{-1})$ has a pole. 
Let
\[ M^{-1}(\lambda) = \frac{1}{\lambda}M_{-1} + M_0 + \lambda M_1 + O(\lambda^2); 
\]
the only way that $\mbox{trace}(M^{-1})$ can fail to have a pole is if
$\mbox{trace}(M_{-1}) = 0$. If this happens then $M_{-1}$ must be of the form
\[ M_{-1} = \left(\begin{array}{cc} \beta & \gamma \\ \gamma & -\beta 
\end{array}\right); \]
recalling that $M_{-1}$ is a real non-zero matrix this implies that $M_{-1}$
is of full rank, and hence $M$ must have a zero rather than a pole when 
$\lambda=0$. This contradiction proves that $\mbox{trace}(M^{-1})$ has 
a pole, and hence that $\Psi$ has a removable singularity.

We now turn to the case $n>2$. We know that $M$ and $M^{-1}$ are
analytic functions of $\lambda$ with singularities at $\lambda=0$.
Let $\mu_1(\lambda),\ldots,\mu_n(\lambda)$ be the eigenvalues of
$\lambda M^{-1}(\lambda)$. As $\lambda M^{-1}(\lambda)$ is analytic, and 
symmetric 
in the sense  of Kato \cite[p. 120]{kn:Kato}, we know from the remark at the 
bottom of 
page 121 in \cite{kn:Kato} that $M^{-1}$ has an analytic Schur decomposition of 
the form
\[ M^{-1}(\lambda)R(\lambda) = R(\lambda) \frac{1}{\lambda}D(\lambda) \]
on a punctured neighbourhood of $\lambda = 0$. Here
\[ D=\mbox{diag}(\mu_1,\ldots,\mu_n), \]
where $\mu_1,\ldots,\mu_n$ are analytic at $\lambda = 0$, while the matrix 
$R(\lambda)$ is 
analytic at $\lambda = 0$ and is real orthogonal ($R^{-1} = R^T$) for all 
sufficiently small 
real $\lambda$. This orthogonality of $R$ for real $\lambda$ means that $R^{-1}$ 
is also 
analytic at $\lambda = 0$. To see this, observe that the only type of 
singularity which 
$R^{-1}$ could have would be a pole. A pole would cause $R^{-1}(\lambda)$ to 
blow up as 
$\lambda$ approached zero through real values, contradicting the regularity of 
$R$ by the 
orthogonality $R^{-1} = R^T$ for real $\lambda$.

The Schur decomposition of $\Psi$ is clearly
\be \Psi(\lambda) = R(\lambda) (\lambda^{-1}D(\lambda) + \alpha I)^{-1} 
R^{-1}(\lambda), 
 \label{eq:mmrev7} \ee
and the eigenvalues of $\Psi$ are clearly
\[ (\alpha + \frac{1}{\lambda}\mu_1)^{-1},\ldots, 
(\alpha+\frac{1}{\lambda}\mu_n)^{-1}. \]
These are all analytic functions of $\lambda$. If $j$ is such
that $\mu_j \neq 0$  at $\lambda = 0$ then the corresponding
eigenvalue of $\Psi$ clearly has a zero at $\lambda=0$. If 
$j$ is such that $\mu_j$ has a zero of order at least 2 at $\lambda=0$, then the
corresponding eigenvalue of $\Psi$ has a removable singularity
at $\lambda=0$ provided $\alpha\neq 0$. If $j$ is such that
$\mu_j$ has a simple zero  at $\lambda=0$ then the
corresponding eigenvalue of $\Psi$ will have a removable singularity
at $\lambda=0$ for all but one value of $\alpha$. In particular, since
$\mu_j$ is real-valued for real $\lambda$, the corresponding eigenvalue of 
$\Psi$ has a removable singularity for $\Im(\alpha)\neq 0$. Whenever the 
eigenvalues of 
$\Psi$ are analytic, so is $\Psi$ itself, from the Schur decomposition 
(\ref{eq:mmrev7}) in which $R(\lambda)$ and $R^{-1}$ are analytic at $\lambda = 
0$.
This completes the proof.
\vspace{2mm}

\noindent {\bf Remark} The proof for $n>2$ can be extended to show
that under the hypothesis of the Bennewitz conjecture -- namely,
that the ranks of the residue matrices of $M$ and $M^{-1}$ sum to
$n$ -- the matrix $\Psi$ has a removable singularity for any non-zero
$\alpha$. In other words, the result of the case $n=2$ is recovered
in this special case. \hfill  $\Box$
\vspace{2mm}

\subsection{The initial value problem for $\Psi$}
\label{section:3.2}
In the rest of this section we shall consider the case $n=2$: that of the
fourth order Sturm-Liouville problem. 
We start by recalling the method proposed for the computation of the matrix
$M$ in \cite{kn:lotsofus}. An interval $[0,X]$ is chosen, with $X$ suitably
large; the fourth order Sturm-Liouville equation is cast in the form
\[ J z' = S z, \]
where $J$ is the symplectic matrix
\[ J = \left(\begin{array}{rrrr} 0 & 0 & -1 & 0 \\
                                 0 & 0 & 0 & -1 \\
                                 1 & 0 & 0 & 0 \\
                                 0 & 1 & 0 & 0 \end{array} \right), \]
$S$ is the symmetric matrix
\be S = \left(\begin{array}{rrrr} \lambda w - q & 0 & 0 & 0 \\
                                 0 & -s & 1 & 0 \\
                                 0 & 1 & 0 & 0 \\
                                 0 & 0 & 0 & 1/p \end{array} \right), 
\label{eq:mm15} \ee
and $z$ is the vector of quasi-derivatives
\[ z = \left( \begin{array}{c} y \\ y' \\ -(py'')'+sy' \\ py'' \end{array} 
\right). \]
Then we consider the matrix initial value problem consisting of the differential
equation
\[ J Z' = S Z \]
and the initial condition
\[ Z(X) = \left(\begin{array}{cc} 0 & 0 \\
                                  0 & 0 \\
                                  1 & 0 \\
                                  0 & 1 \end{array} \right). \]
Let $Z(x)$ denote the solution of this problem, where $0 \leq x \leq X$.
We partition $Z$ as
\[ Z(x) = \left(\begin{array}{c} U(x) \\
                                 V(x) \end{array} \right), \]
and form the corresponding initial value problem for the matrix $UV^{-1}$.
We solve this initial value problem, starting from $x=X$, to find $UV^{-1}(0)$.
Our approximation to $M$ is then given by
\be M \approx (UV^{-1}(0))^{-1}. \label{eq:mm14} \ee
(The formula would be exact if we had $X=+\infty$.)
If we replace $M$ in (\ref{eq:mm1}) by the expression on the right hand
side of (\ref{eq:mm14}) then we get
\[ \Psi \approx (\alpha I + UV^{-1}(0))^{-1}. \]
Clearly, then, the process of approximating $\Psi$ can be reduced to that of
deriving an initial value problem for the matrix
\[ \Gamma(x) := (\alpha I + UV^{-1}(x))^{-1}. \]
The initial condition is obvious: $\Gamma(X) = \alpha^{-1}I$. The differential
equation is also quite straightforward to derive. If the matrix $S$ in 
(\ref{eq:mm15}) is partitioned as 
\[ S = \left(\begin{array}{cc} S_{1,1} & S_{1,2} \\ S_{2,1} & S_{2,2} 
\end{array}
 \right), \]
then it turns out that
\be \Gamma' = -\left\{ \Gamma S_{2,1} (I-\alpha\Gamma) + 
(I-\alpha\Gamma)S_{1,2}\Gamma
 + (I-\alpha\Gamma)S_{1,1}(I-\alpha\Gamma) + \Gamma S_{2,2}\Gamma \right\}.
 \label{eq:mmrev3} \ee
We solve this equation using the NAG routine D02QGF, which allows reverse
communication for evaluation of the right hand side of the differential 
equation:
this helps to keep the programme structure simple when the right hand side is
complicated. D02QGF is a variable-order, variable-step Adams code and is 
therefore 
able to cope with mild stiffness. In practice, we noted in \cite{kn:lotsofus}
that stiffness is not usually a problem unless $X$ has been chosen much
larger than necessary.

\subsection{Determining the Taylor expansion of $\Psi$}
Determining an approximate Taylor expansion of an analytic function from 
numerical
values of the function is not easy. The number of coefficients in the
expansion which can be computed reliably depends on the accuracy with
which the function values can be computed, on the rate of decay of the
Taylor coefficients as one proceeds up the series and, ultimately, on the
precision of the machine arithmetic. 

Our problem is slightly compounded by the fact that we have a function $\Psi$
with a removable singularity at the point around which we wish to
expand it ($\lambda=0$). We cannot compute $\Psi(0)$; indeed we cannot
compute $\Psi(\lambda)$ for any $\lambda$ with zero imaginary part.
Our approach has been to compute $\Psi$ at a sequence of points 
\[ \lambda = \frac{\mu}{2^{j}}, \;\;\; j=0,1,2,\ldots, \]
where $\mu$ is a fixed complex number with $\Im(\mu)\neq 0$, and solve
a Vandermonde system (using the algorithm of Bj\"{o}rck and Pereyra 
\cite{kn:Bjorck})
to obtain approximations to the Taylor coefficients.

We shall now consider how the different sources of error and the 
ill-conditioning
of the Vandermonde system will affect the approximations to the Taylor 
coefficients 
which we obtain. For simplicity we shall set aside our matrix-valued function
$\Psi$ and consider a complex-valued function $f$ given by a Taylor
expansion
\be f(z) = a_0 + a_1 z + a_2 z^2 + \cdots + a_m z^m + a_{m+1}z^{m+1} + 
  \cdots  \label{eq:mm10} 
\ee
Clearly the following system of equations holds:
\be
\begin{array}{rcl}
 a_{0} + a_{1}\mu + a_{2}\mu^2 + \cdots a_{m}\mu^{m} & = & g(\mu) \\
 a_{0} + a_{1}\frac{\mu}{2} + a_{2}\left(\frac{\mu}{2}\right)^2 + 
 \cdots a_{m}\left(\frac{\mu}{2}\right)^{m} & = & g(\mu/2) \\
 a_{0} + a_{1}\frac{\mu}{4} + a_{2}\left(\frac{\mu}{4}\right)^2 + 
 \cdots a_{m}\left(\frac{\mu}{4}\right)^{m} & = & g(\mu/4) \\
 \cdots & & \cdots \\
 a_{0} + a_{1}\frac{\mu}{2^n} + a_{2}\left(\frac{\mu}{2^m}\right)^2 + 
 \cdots a_{m}\left(\frac{\mu}{2^m}\right)^{m} & = & g(\mu/2^m) 
\end{array} \label{eq:mm11}
\ee
where the function $g$ is given by
\[ g(z) = f(z) - a_{m+1}z^{m+1} - a_{m+2}z^{m+2} - \cdots. \]
A Vandermonde matrix is an $(m+1)\times (m+1)$ matrix of the form
\[ V = \left(\begin{array}{ccccc} 1 & 1 & 1 & \cdots & 1 \\
                           \alpha_0 & \alpha_1 & \alpha_2 & \cdots & \alpha_m \\
                           \alpha_0^2 & \alpha_1^2 & \alpha_2^2 & \cdots & 
\alpha_m^2 \\
                           \cdot & \cdot & \cdot & \cdots & \cdot \\
                           \cdot & \cdot & \cdot & \cdots & \cdot \\
                           \alpha_0^m & \alpha_1^m & \alpha_2^m & \cdots & 
\alpha_m^m 
             \end{array}\right),
\] 
where $\alpha_0,\ldots,\alpha_m$ are distinct complex numbers. If we 
let ${\bf a} = (a_0,a_1\mu,a_2\mu^2,\ldots,a_m\mu^m)^T$ and
${\bf g} = (g(\mu/2^m),g(\mu/2^{m-1}),\ldots,g(\mu))^T$ then we can cast our 
system as a
dual Vandermonde problem of the form
\[ V^T {\bf a } = {\bf g}, \]
where $\alpha_0 = 1/2^m$, $\alpha_1 = 1/2^{m-1}$,\ldots,$\alpha_m = 1$. Since 
the
$\alpha_j$ are positive real numbers arranged in ascending order, the error 
analysis of
Higham \cite{kn:Higham} is now applicable to the solution of this system by the
Bj\"{o}rck-Pereyra algorithm. In particular this error analysis shows, when we 
know
${\bf g}$ `exactly', that the Bj\"{o}rck-Pereyra algorithm  introduces 
essentially no more error 
into ${\bf a}$ than is already implied by the storage of ${\bf g}$ in machine 
arithmetic. Since we do not have any of the special sign properties on the 
vector
${\bf g}$ which would make for a better error estimate, this suggests that the
contribution of the machine precision $\epsilon$ to the error (measured in the 
norm
$\| \cdot \|_1$) will be a term of the order
\be m\epsilon \| V^{-1} \|_{\infty} \leq C m \epsilon 2^{m^2}, \label{eq:NH} \ee
where $C$ is independent of $m$ \cite{kn:Higham}.
With a machine precision of $10^{-16}$ this suggests that, on grounds of 
roundoff
alone, it will not be possible to obtain reasonable accuracy in the vector
${\bf a}$ for $m$ much greater than 6. This was borne out in the experiments
which we conducted. Of course we could use $\| \cdot \|_\infty$ to measure the 
error, instead of $\| \cdot \|_1$; however this would make no difference, as the 
ratio $\| {\mathbf x} \|_1/\| {\mathbf x} \|_\infty$ is never greater than $m$ 
for any non-zero $m$-vector ${\bf x}$, and we have already seen that $m$ must 
be quite small.

We now turn to the contribution to the error arising from ${\bf g}$: we 
do not know ${\bf g}$ exactly because we do not know $g$ exactly.
We must approximate $g$ by $f$. This entails an error
\be g(z) - f(z) = - a_{m+1}z^{m+1} - a_{m+2}z^{m+2} - \cdots, \label{eq:mmrhs} 
\ee
plus a further term arising from the numerical integration of the differential
equation satisfied by $f$.
%
% If the  error arising from the numerical integration of
%the differential equation had no special structure, its effect on the whole 
process
%might be similar to contaminating $g$ with a roundoff error of order $\tau$, 
where
%$\tau$ would (optimistically) be of the order of the local truncation error of 
the
%integrator used. Thus by analogy with (\ref{eq:NH}) we might expect the error
%arising in the vector ${\bf a}$ as a result of numerical integration with 
accuracy
%$\tau$ to be bounded by a term of the form
%\[ C m \tau 2^{m^2}, \]
%With a machine precision of $10^{-16}$ we would be fortunate to get
%$\tau \approx 10^{-12}$, even using a very small stepsize; this would impose an 
%even stricter bound on $m$ than is implied by the machine arithmetic. This does 
not 
%seem to happen in practice: we seem to see an integration-induced error which 
is much 
%closer to that caused by the right hand side of (\ref{eq:mmrhs}), suggesting 
that the error
%due to the integration is either much smaller than the right hand side of
%(\ref{eq:mmrhs}) or else has some special structure such as the sign 
alternation 
%proposed by Higham \cite{kn:Higham}. 
Neglecting the integration error for the moment, we observe that 
%
%Finally, then, we are able to turn to the error in the coefficients 
$a_0$,\ldots,$a_n$ 
%caused purely by the replacement of $g$ by f$. To do this we observe that, 
%
from (\ref{eq:mm11}), the coefficients $a_0$,\ldots,$a_m$ are the interpolation 
coefficients
for the function $g$ at the points $\mu$, $\mu/2$,\ldots,$\mu 2^{-m}$. As 
interpolation
is a linear process, 
\begin{eqnarray}
\mbox{Interpolant of $g$} & = & \mbox{Interpolant of $f$} + \mbox{Interpolant of 
$(g-f)$} 
 \nonumber \\
& = & \mbox{Interpolant of $f$} - \sum_{p=1}^{\infty}a_{m+p}\{\mbox{Interpolant 
of $z^{m+p}$}\},
 \label{eq:mmrev5} 
\end{eqnarray}
the second equality coming from (\ref{eq:mmrhs}). Now let
\[ \sum_{\nu=0}^{m}c_{\nu}^{(p)}\left(\frac{z}{\mu}\right)^\nu \]
denote the $m$th degree polynomial interpolant to the function $z\mapsto 
z^{m+p}$ at the 
points $1$, $2^{-1}$, \ldots, $2^{-m}$. By (\ref{eq:mmrev5}) the change 
in the computed value of $a_{\nu}\mu^{\nu}$ caused by approximating $g$ by $f$ 
is given by
\[ \sum_{p=1}^{\infty}a_{m+p}c_{\nu}^{(p)}. \]
The coefficients $c_{\nu}^{(p)}$ can be computed explicitly using the
Bj\"{o}rck-Pereyra algorithm. The solution is given by
\be c_{\nu}^{(p)} = 
\sum_{j=0}^{m-\nu}(-1)^{m-\nu-j}\beta_{m-\nu-j}(1,2^{-1},\ldots,2^{-m+1})
 \mu^{m+p}\left(\frac{1}{2^j}\right)^{m+p-j}
 \frac{\prod_{r=0}^{j-1}(1-2^{m+p-r})}{\prod_{r=1}^{j}(1-2^r)},
\label{eq:mmcp} \ee
where $\beta_k(y_0,\ldots,y_{m-1})$ denotes the sum of all products of
$k$ distinct elements of the set $\{ y_0,\ldots, y_{m-1} \}$. We use this
formula to get a bound on $c_{\nu}^{(p)}$. First, we consider
$\beta_{k}(1,2^{-1},\ldots,2^{-m+1})$. This is a sum of $\left( \begin{array}{c} 
m \\ k
 \end{array}\right)$ terms, of which the greatest is $1.2^{-1}\ldots 2^{-k+1} = 
2^{-k(k-1)/2}$
Thus
\be \beta_{m-\nu-j}(1,2^{-1},\ldots,2^{-m+1}) \leq \left( \begin{array}{c} m \\ 
m - \nu -j 
 \end{array}\right) 2^{-(m-\nu-j)(m-\nu-j-1)/2}. \label{eq:mmcp1} \ee
Next we tackle the term
\begin{eqnarray} \left| \prod_{r=0}^{j-1}(1-2^{m+p-r}) \right | & = & 
 \prod_{r=0}^{j-1}2^{m+p-r}(1-2^{-m-p+r}) \nonumber \\
 & \leq & \prod_{r=0}^{j-1}2^{m+p-r} \nonumber \\
 & = & 2^{(m+p-(j-1)/2)j} \label{eq:mmcp2} 
\end{eqnarray}
Similarly,
\begin{eqnarray}
\left| \prod_{r=1}^{j}(1-2^r) \right| & = & \prod_{r=1}^{j}2^r(1-2^{-r}) 
\nonumber \\
 & = & 2^{j(j+1)/2}\prod_{r=1}^{j}(1-2^{-r}) \nonumber \\
 & \geq & 2^{j(j+1)/2}\prod_{r=1}^{j}\exp(-2.2^{-r}) \nonumber \\
 & \geq &  2^{j(j+1)/2}\exp(-2). \label{eq:mmcp3} 
\end{eqnarray}
Substituting (\ref{eq:mmcp1}), (\ref{eq:mmcp2}) and (\ref{eq:mmcp3}) back into
(\ref{eq:mmcp}) we obtain the estimate
\be
 |c_{\nu}^{(p)}| \leq \sum_{j=0}^{m-\nu}\left(\begin{array}{c} m \\ m - \nu - j 
\end{array}\right)
 2^{-(m-\nu-j)(m-\nu-j-1)/2}|\mu|^{m+p} 
\left(\frac{1}{2^j}\right)^{m+p-j}\exp(2)\frac{
 2^{(m+p-(j-1)/2)j}}{2^{j(j+1)/2} } 
\ee
which simplifies to give
\be |c_{\nu}^{(p)}| \leq \sum_{j=0}^{m-\nu}\left(\begin{array}{c} m \\ m - \nu - 
j \end{array}\right)
 |\mu|^{m+p}\exp(2)2^{-(m-\nu)^2/2+(m-\nu)(2j+1)/2-j(j+1)/2}.
\ee
By the change of dummy summation variable $j = m - \nu - k$, this gives
\be |c_{\nu}^{(p)}| \leq \sum_{k=0}^{m-\nu}\left(\begin{array}{c} m \\ k 
\end{array}\right)
 |\mu|^{m+p}\exp(2)2^{-k(k-1)/2}.
\ee
For our purposes, it suffices to make a very blunt estimate at this stage: throw
away the powers of 2, and extend the summation up to $k=m$. Since
\[ \sum_{k=0}^{m}\left(\begin{array}{c} m \\ k \end{array}\right) = 2^m, \]
we get the bound
\[ |c_{\nu}^{(p)}| \leq |\mu|^{p}|2\mu|^{m}\exp(2) \leq \exp(2)|2\mu|^{m+p}. \]
Thus, neglecting the integration error, the change in the computed value of 
$a_{\nu}\mu^{\nu}$ 
caused by replacing $g$ by $f$ will be bounded by
\be \exp(2)\sum_{p=1}^{\infty}|a_{m+p}||2\mu|^{m+p}. \label{eq:mmpert} \ee
Provided $2\mu$ is strictly within the radius of convergence of the power 
series,
this will tend to zero as $m$ tends to infinity; moreover, if $a_{m+1}\neq 0$,
the leading order term for small $\mu$ will be
\[ |a_{m+1}|\exp(2)|2\mu|^{m+1}. \]
In order to assess the likely order of magnitude for a suitable value of $\mu$ 
it
seems reasonable to ask that this error term be of the same order of magnitude
as the error arising from roundoff. The constant $C$ in (\ref{eq:NH}) is 
independent
of $m$ so we can neglect it; the coefficient $a_{m+1}$ we obviously do not know,
but if we assume that it is $O(1)$ then we obtain 
\[ m. 2^{m^2}\epsilon \approx |2\mu|^{m+1}. \]
We have already seen that with $\epsilon = 10^{-16}$, the choice of $m=6$ is 
likely
to be the biggest possible. In some sense, this modest value of $m$ justifies 
the
assumption that $a_{m+1}$ is an $O(1)$ quantity. It also gives 
$m. 2^{m^2}\epsilon \approx 4\times 10^{-5}$,
which suggests $\mu \approx 0.1$. With $m=5$, on the other hand, we have
$m.2^{m^2}\epsilon \approx 1.7\times 10^{-8}$, which suggests $\mu\approx 
0.025$.

Based on these observations we devised
the following algorithm for computing the first $k+1$ Taylor coefficients
of the matrix $\Psi$, where $k < 5$.
\begin{description}
\item[1.] Make the tolerance $TOL$ for the computation of $\Psi$ as small as 
possible 
 within the constraints of reasonable run-times. This depends on the machine at
 one's disposal. 
\item[2.] Start with $|\mu|\approx 0.025$ (say) and $n=k+1$.
\item[3.] Compute approximations to $a_0,\ldots,a_{m}$; from these extract 
 approximations to $a_0,\ldots,a_k$, the coefficients of interest.
\item[4.] If $m < 7$, increase the value of $m$ and compute new coefficient 
 approximations. 
\item[5.] While the approximations seem to be improving and $m<7$, keep 
 increasing $m$ and computing new approximations.
\item[6.] When the sequence of coefficient approximations appears to start
 to diverge stop increasing $m$. Discard the latest (starting-to-diverge) 
 approximations.   
\item[7.] Now regard $m$ as fixed and start to double $\mu$. Follow the
 same process as above, doubling $\mu$ while this seems to improve the
 values of $a_0,\ldots,a_k$. Stop either when the user's target accuracy
 is achieved, or when the approximations seem to start to diverge,
 or when a doubling of $\mu$ would give $|\mu| > 0.5$. Return a warning
 flag ($IFAIL = 2$) if the target accuracy has not been reached.
\end{description}

The error due to integration is never explicitly controlled in this
process, though step {\bf 6} should ensure that $m$ is never taken large
enough to magnify the integration error to an unacceptable level.

Typically, for computing the Taylor coefficients of our matrix $\Psi$, we might 
start with $\mu = 0.025i$, and compute $\Psi$ using an initial value solver with 
$TOL=10^{-11}$. If $\Psi$ has a Taylor expansion
\[ \Psi(\lambda) = \Psi_{0} + \lambda\Psi_{1} + \lambda^2\Psi_{2} 
 + \lambda^3\Psi_{3} + \lambda^4\Psi_{4} + \cdots, \]
we usually find that only the coefficients 
$\Psi_{0},\ldots,\Psi_{3}$ can be computed with an accuracy of $10^{-4}$
or better; the accuracy of $\Psi_3$ might be $10^{-4}$, of $\Psi_{2}$
about $10^{-6}$, of $\Psi_{1}$ about $10^{-8}$, while $\Psi_{0}$ might
have an accuracy of $10^{-10}$, achieved with $m=6$ or $m=7$ and with
a value of $\mu$ of about $0.2i$. The precise details depend on the problem
in question. The deterioration in the accuracy of the coefficients as one 
proceeds up the series need not be a problem if one intends to use them 
simply to compute values of $\Psi(\lambda)$ by a truncated Taylor expansion for
small values of $|\lambda|$.

\subsection{Recovering the residue of $M$ from $\Psi$}
We suppose that the first few terms of the Taylor expansion of $M$ have
been determined:
\be  \Psi(\lambda) = \Psi_{0} + \lambda\Psi_{1} + \lambda^2\Psi_{2} + 
                    \lambda^3\Psi_3 + \cdots \label{eq:mm3}
\ee
We want to determine the first few coefficients in the Laurent expansion
of $M$:
\be M(\lambda) = \frac{1}{\lambda}M_{-1} + M_{0} + \lambda M_{1} + \lambda^2 
M_{2}
 + \cdots \label{eq:mm4} 
\ee
As we shall see, from the first $n$ terms in the Taylor expansion of 
$\Psi$ we can determine at most the first $n-1$ terms in the Laurent expansion
of $M$. Thus, although we have eliminated the problems associated with
trying to compute $M$ near the pole, we have paid a price in terms of having
to compute more Taylor coefficients than we get repaid in Laurent coefficients.

Equation (\ref{eq:mm1}) may be rearranged to yield
\begin{eqnarray*} M^{-1} & = & \Psi^{-1} - \alpha I \\
                         & = & \frac{\Psi^{A}}{\det(\Psi)} - \alpha I \\
                         & = & \frac{\Psi^{A}-\alpha\det(\Psi)I}{\det(\Psi)}.
\end{eqnarray*}
Inverting both sides,
\begin{eqnarray*} M & = & \det(\Psi) (\Psi^{A}-\alpha\det(\Psi)I)^{-1} \\
           & = & \det(\Psi) \frac{\Psi - \alpha \det(\Psi) I}{
                              \det(\Psi - \alpha \det(\Psi) I)}
\end{eqnarray*}
Now we expand the determinant in the denominator to get
\[ \det(\Psi - \alpha \det(\Psi) I) = 
\det(\Psi)\left\{1-\alpha\mbox{trace}(\Psi)
 + \alpha^2\det(\Psi)\right\}. 
\]
Thus we obtain the formula
\be M = \frac{\Psi - \alpha \det(\Psi) I}{1 - \alpha \mbox{trace}(\Psi) +
 \alpha^2 \det(\Psi)}. \label{eq:mm5}
\ee
We shall obtain the Laurent series for $M$ by Taylor expansion of the
numerator and denominator in (\ref{eq:mm5}). From (\ref{eq:mm3}) we obtain
the expansions
\be \mbox{trace}(\Psi) = \mbox{trace}(\Psi_{0}) + \lambda \mbox{trace}(\Psi_{1})
 + \lambda^2 \mbox{trace}(\Psi_{2}) + \lambda^3 \mbox{trace}(\Psi_{3}) +
 \cdots \label{eq:mm6}, 
\ee
\begin{eqnarray*} \det(\Psi) & = & \det\Psi_{0} + 
\lambda\mbox{trace}(\Psi_{0}^{A}\Psi_{1}) \\
 &  + & \lambda^2 (\det(\Psi_{1}) + \mbox{trace}(\Psi_{0}^{A}\Psi_{2}) ) \\
 & + & \lambda^3 (\mbox{trace}(\Psi_{1}^{A}\Psi_{2}) + 
\mbox{trace}(\Psi_{0}^{A}\Psi_{3}) ) \\
 & + & \lambda^4 (\mbox{trace}(\Psi_{0}^{A}\Psi_{4}) + 
\mbox{trace}(\Psi_{1}^{A}\Psi_{3}) 
              + \det(\Psi_{2})) \\
 & + &  \cdots
\end{eqnarray*}
We know that the denominator in (\ref{eq:mm5}) must have a zero at $\lambda=0$ 
because
$M$ has a pole at $\lambda = 0$ by hypothesis: thus
\begin{eqnarray} 1 - \alpha \mbox{trace}(\Psi) + \alpha^2\det(\Psi) & = &
 \lambda(\alpha^2 \mbox{trace}(\Psi_{0}^{A}\Psi_{1}) - \alpha 
\mbox{trace}(\Psi_{1})) \nonumber \\ 
 & + & \lambda^2(\alpha^2\det(\Psi_{1}) + 
\alpha^2\mbox{trace}(\Psi_{0}^{A}\Psi_{2})
 - \alpha \mbox{trace}(\Psi_{2})) \nonumber \\
 & + & 
\lambda^3(\alpha^2\mbox{trace}(\Psi_{0}^{A}\Psi_{3})+\alpha^2\mbox{trace}(
        \Psi_{1}^{A}\Psi_{2}) - \alpha\mbox{trace}(\Psi_{3})) \nonumber \\
 & + & 
\lambda^4(\alpha^2\mbox{trace}(\Psi_{0}^{A}\Psi_{4})+\alpha^2\mbox{trace}(
        \Psi_{1}^{A}\Psi_{3}) +\alpha^2\det(\Psi_{2})- 
\alpha\mbox{trace}(\Psi_{4}))
\nonumber \\
 & + & \cdots \label{eq:mm7a} 
\end{eqnarray}
For brevity we shall write this expression in the form 
\be 
 1 - \alpha \mbox{trace}(\Psi) + \alpha^2\det(\Psi) = \lambda a_{1} + 
 \lambda^2 a_{2} + \lambda^3 a_{3} + \lambda^4 a_{4} + \cdots, \label{eq:mm7} 
\ee
where the coefficients $a_{1}$, $a_{2}$, $a_{3}$ and $a_{4}$ are evident by 
comparing
(\ref{eq:mm7a}) and (\ref{eq:mm7}). We shall also write
\be \Psi - \alpha\det(\Psi)I  =  A_{0} + \lambda A_1 + \lambda^2 A_2 + \lambda^3 
A_3
 + \lambda^4 A_4 + \cdots, \label{eq:mm8}
\ee
where 
\begin{eqnarray*}
 A_{0} & = & \Psi_{0} - \alpha\det(\Psi_{0})I, \\
 A_{1} & = & \Psi_{1} - \alpha\mbox{trace}(\Psi_{0}^{A}\Psi_{1})I, \\
 A_{2} & = & \Psi_{2} - \alpha\left(\det(\Psi_{1})+
              \mbox{trace}(\Psi_{0}^{A}\Psi_{2})\right)I, \\
 A_{3} & = & \Psi_{3} - \alpha\left(\mbox{trace}(\Psi_{0}^{A}\Psi_{3})+
              \mbox{trace}(\Psi_{1}^{A}\Psi_{2})\right)I, \\
 A_{4} & = & \Psi_{4} - \alpha\left(\mbox{trace}(\Psi_{0}^{A}\Psi_{4})+
              \mbox{trace}(\Psi_{1}^{A}\Psi_{3})+\det(\Psi_{2})\right)I. \\
\end{eqnarray*} 
If $a_1\neq 0$, then we may combine (\ref{eq:mm5}), (\ref{eq:mm7}) and 
(\ref{eq:mm8}) to gives us the Laurent expansion
\begin{eqnarray} M(\lambda) & = & \lambda^{-1}\frac{A_{0}}{a_{1}} + 
\left(\frac{A_{1}}{a_{1}}-
 \frac{a_{2}}{a_{1}^{2}}A_{0}\right) \nonumber \\
 & + & \lambda \left(\frac{A_{2}}{a_{1}} - \frac{a_{2}}{a_{1}^{2}}A_{1} +    
\left(\frac{a_{2}^{2}}{a_{1}^{3}}-\frac{a_{3}}{a_{1}^{2}}\right)A_{0}\right) 
\nonumber \\
 & + & \lambda^2 \left(\frac{A_{3}}{a_{1}} - \frac{a_{2}}{a_{1}^{2}}A_{2} +
       \left(\frac{a_{2}^{2}}{a_{1}^{3}}-\frac{a_{3}}{a_{1}^{2}}\right)A_{1}+       \left(\frac{(2a_{1}a_{3}-a_{2}^{2})a_{2}}{a_{1}^{4}}-\frac{a_{4}}{a_{1}^{2}}\right)A_{0}
       \right) \label{eq:mm9} \\
 & + & \cdots \nonumber
\end{eqnarray}
Since $M$ has at worst a simple pole, when $a_1=0$ then $A_0=0$. In this 
case the Laurent expansion of $M$ becomes
\begin{eqnarray} M(\lambda) & = & \lambda^{-1}\frac{A_{1}}{a_{2}} + 
 \left(\frac{A_{2}}{a_{2}}- \frac{a_{3}}{a_{2}^{2}}A_{1}\right) \nonumber \\
 & + & \lambda \left(\frac{A_{3}}{a_{2}} - \frac{a_{3}}{a_{2}^{2}}A_{2} +       
\left(\frac{a_{3}^{2}}{a_{2}^{3}}-\frac{a_{4}}{a_{2}^{2}}\right)A_{1}\right) 
 \label{eq:mm9a} \\
 & + & \cdots \nonumber
\end{eqnarray}
This happens when $\Psi(0) = \frac{1}{\alpha}I$, which happens when the
residue matrix $\mbox{Res}(M,0)$ has full rank, giving $M^{-1}$ a 
zero at $\lambda=0$. The code checks that the value of $a_1$, and uses
(\ref{eq:mm9}) if $|a_1| > TOL$, (\ref{eq:mm9a}) if $|a_1| < TOL$, where
$TOL$ is the tolerance used by $D02QGF$ in the computation of $\Psi$. The
case $a_1 = a_2 = 0$ cannot arise with $2\times 2$ matrices: for $a_2=0$ 
necessarily 
implies $A_1=0$. Since $\Psi(0) = \frac{1}{\alpha}I$, we see that $\det \Psi(0) 
\neq 0$;
this implies that $M^{-1} = (\Psi^{A}-\det(\Psi)I)/\det(\Psi)$ has a double zero 
at 
$\lambda=0$, implying that $M$ has a double pole (or worse). This is not 
possible, 
as both $M$ and $M^{-1}$ have, at worst, simple 
poles.

\section{Numerical Experiments}
\label{section:numerics}
Our primary objective in these experiments was to compute the residues of
Titchmarsh-Weyl matrices for a number of fourth order Sturm-Liouville 
equations and to use these, together with the Bennewitz Conjecture, to 
decide whether HELP inequalities hold for these equations.

Before listing our example problems, we mention the following
useful result. To avoid complicated conditions on quasiderivatives we state
the result for smooth coefficients in the differential operator.
\begin{lemma}\label{lemma:new}
Suppose that ${\cal L}$ is a fourth order differential operator on the domain of 
functions $f$ in $L^2[0,\infty)$ which are four times continuously 
differentiable
and are such that ${\cal L}f \in L^2[0,\infty)$. Suppose that ${\cal L}$ has
the form ${\cal L} = \ell^2$ where $\ell$ is a second order operator
\[ \ell(y)(x) = -y''(x) + q(x) y(x) \]
and $q$ is twice continuously differentiable with $q(0) = 0$, $q'(0) = 0$.
Suppose also that ${\cal L}$ is strong limit-point at infinity.
Then any eigenfunction $y$ of ${\cal L}$ subject to Neumann boundary
conditions $y'''(0) = 0 = y''(0)$ which is not in the null-space
of $\ell$ generates an eigenfunction $z=\ell y$ of ${\cal L}$
subject to Dirichlet conditions $z(0) = 0 = z'(0)$.
\end{lemma}

\noindent {\bf Proof}\, Suppose $y$ is as described, so that
${\cal L}y = \lambda y$ for some real $\lambda$, and let $z=\ell y$.
Because $y$ is not in the null-space of $\ell$, $z$ is non-trivial.
Also $z(0)$ $=$ $-y''(0)$ $+$ $q(0) y(0)$ $ = $ $ 0 $ because
$y''(0) = 0$ and $q(0) = 0$, and $z'(0) $ $ = $ $ -y'''(0) $ $ + $
$ q(0) y'(0) $ $+$ $q'(0) y(0)$ $=$ $0$ since $y'''(0) = 0$
and $q(0)$ $=$ $q'(0)$ $=$ $0$. Thus $z$ satisfies the Dirichlet
boundary conditions. Clearly ${\cal L}z$ $=$ $\ell^2 \ell y$ $=$
$\ell {\cal L}y$ $=$ $\ell \lambda y$ $=$ $\lambda \ell y$ $=$ 
$\lambda z$, so $z$ satisfies the differential equation
${\cal L}z$ $=$ $\lambda z$: as the coefficient $q$ is twice
continuously differentiable, this makes it easy to see that $z$
is four times continuously differentiable. Finally, $z$ is square integrable.
This follows because
\[ \langle z,z \rangle = \langle \ell y, \ell y \rangle 
 = \langle \ell^2 y, y \rangle = \lambda \langle y,y \rangle, \]
the penultimate equality using the fact that $\ell^2 = {\cal L}$ is 
strong limit-point at infinity together with the fact that 
$\ell y(0)$ $=$ $(\ell y)'(0)$ $=$ $0$. \hfill $\Box$

\begin{description}

\item[Equation 1] The differential equation
\[ y^{(iv)} - (s(x)y')' + q(x) y = \lambda y \]
on the interval $[0,\infty)$, with coefficients 
\[ s(x) = \frac{8x^2(x^4-3x^2-5)}{(x^2+1)^2}, \]
\[ q(x) = \frac{4[4x^{12}-24x^{10}-7x^8+96x^6+46x^4-60x^2-15]}{(x^2+1)^4}. \]
This equation is strong limit-point at infinity. Imposing a Dirichlet
boundary condition $y(0) = y'(0) = 0$ reveals that this problem was carefully
crafted so that $\lambda = 0$ would be an eigenvalue of multiplicity 2; the
reader may check that
\[ y(x) = (x^2+x^4)\mbox{e}^{-x^2}, \;\;\; y(x) = (x^3+x^5)\mbox{e}^{-x^2} \]
are the corresponding eigenfunctions. We arranged this because we suspected 
that it would result in a problem for which we would have
\be \mbox{rank}(M_D) = 2, \;\;\; \mbox{rank}(M_N) = 0, \label{eq:na1} \ee
although in fact we do not know of any result which would guarantee this.
The numerical results in Table \ref{table:1} suggest that (\ref{eq:na1}) is 
indeed true (the determinant of $M_D$ vanishes nowhere). If it is, 
then the hypotheses of the Bennewitz conjecture are 
satisfied and we have an equation for which a HELP inequality holds.
We should mention that there seems to be a dearth of fourth order
examples with multiple eigenvalues in the literature: indeed, we could
not find any.

Notice that we carried out the computations for two different values
of $\alpha$, to provide an additional check on our results; we also
quote what the code thinks is the imaginary part of the residue matrix.
This ought to be zero, so it provides an indication of the error. We
also quote the error indicator returned by the code: this is reassuringly
of the same order of magnitude as the imaginary part of the computed 
residue matrix.

\begin{table} 
\hspace{0.6in}
\begin{tabular}{|c|} \hline
 Using $\alpha = 1$, truncating $[0,\infty)$ to $[0,100]$: \\ \hline
 $\mbox{Res}(M_D,\lambda=0) = 
   \left(\begin{array}{cc} -384.83 & -167.96 \\
                           -167.96 & -79.29 \end{array}\right) + 
   i \left(\begin{array}{cc} 4.2\times 10^{-2} & -1.7\times 10^{-2} \\
                           -1.7\times 10^{-2} & 6.6\times 10^{-3} 
 \end{array}\right)$ \\ 
 Code error estimate (sup norm): $3.4 \times 10^{-2}$ \\ \hline
 Determinant of residue matrix: $2303.17 + i0.0008$ \\ \hline
 Value of $|a_1|$: $6.3\times 10^{-10}$. \\
 Integration tolerance for $\Psi$: $10^{-9}$ \\ \hline \hline
 Using $\alpha = 1+i$, truncating $[0,\infty)$ to $[0,10]$: \\ \hline
 $\mbox{Res}(M_D,\lambda=0) = 
   \left(\begin{array}{cc} -384.79 & -167.94 \\
                           -167.94 & -79.28 \end{array}\right) + 
   i \left(\begin{array}{cc} 2.5\times 10^{-3} & 1.0\times 10^{-3} \\
                           1.0\times 10^{-2} & 4.1\times 10^{-4} 
 \end{array}\right)$ \\ 
 Code error estimate (sup norm): $2.1 \times 10^{-3}$ \\ \hline
 Determinant of residue matrix: $2303.18 - i0.02$ \\ \hline
 Value of $|a_1|$: $-1.7\times 10^{-12}$. \\
 Integration tolerance for $\Psi$: $10^{-9}$ \\ \hline \hline
\end{tabular} 
\caption{Results for Equation 1}
\label{table:1}
\end{table}

\item[Equation 2] The differential equation
\[ y^{(iv)} - 2(x^2y')'+ (x^4-2)y = \lambda y\]
on the interval $[0,\infty)$. This equation is strong limit-point
at infinity; its differential operator is the formal square of the 
second order operator $\ell y = -y'' + x^2 y$. By Lemma \ref{lemma:new}
we know that all but at most one of the Neumann eigenvalues
will be Dirichlet eigenvalues: in fact if Dirichlet boundary 
conditions $y(0) = 0 = y'(0)$ are imposed, then the eigenvalues are 
$\lambda_k = 16(k+1)^2$, while if Neumann boundary conditions 
$y'''(0)=0=y''(0)$ are imposed, then the eigenvalues are 
$\lambda_k = 16k^2$. Thus at each of the points 
$\lambda = 16(k+1)^2$, $k=0,1,\ldots$, both the Titchmarsh-Weyl 
matrices $M_N$ and $M_D$ will have poles. Given that $M_N = -M_D^{-1}$
it is clear that this means that at each of these poles, the residue of 
$M_N$ and the residue of $M_D$ will both be of rank 1. Together with the 
now-proved Bennewitz Conjecture this means that there are HELP inequalities
(\ref{eq:1.2}) associated with this differential equation, provided the 
operator ${\cal M}$ is defined by
\[ {\cal M}y = y^{(iv)} - 2(x^2y')'+ (x^4-2)y - 16(k+1)^2y, \]
where $k$ is some non-negative integer. This is a result which
Diaz conjectured in his thesis \cite{kn:Dias} but was unable to
prove.

Consulting the numerical results in Table \ref{table:2} we see that our 
code obtains approximations to the residue matrices which are as near
to rank 1 as one could expect: they are matrices whose elements
are not small but whose determinants are  $O(10^{-13})$. 

\begin{table} 
\hspace{0.5in}
\begin{tabular}{|c|} \hline
 Using $\alpha = 1$, truncating $[0,\infty)$ to $[0,20]$: \\ \hline
 $\mbox{Res}(M_D,\lambda=16) = 
   \left(\begin{array}{cc} -82.62549 & -40.74366 \\
                           -40.74366  & -20.09121  \end{array}\right) + 
   i \left(\begin{array}{cc} -7.3\times 10^{-6} & -3.6\times 10^{-6} \\
                           -3.6\times 10^{-6} & -1.8\times 10^{-6} 
 \end{array}\right)$ \\ 
 Code error estimate (sup norm): $3.5 \times 10^{-6}$ \\ \hline
 Determinant of residue matrix: $-4.5\times 10^{-13}-i2.3\times 10^{-18}$  
 \\ \hline
 Value of $|a_1|$: $9.7\times 10^{-3}$. \\
 Integration tolerance for $\Psi$: $10^{-9}$ \\ \hline 
 $\mbox{Res}(M_N,\lambda=16) = 
   \left(\begin{array}{cc} -82.62548 & -40.74366 \\
                           -40.74366  & -20.09121  \end{array}\right) + 
   i \left(\begin{array}{cc} -1.8\times 10^{-5} & -8.7\times 10^{-6} \\
                           -8.7\times 10^{-6} & -4.3\times 10^{-6} 
 \end{array}\right)$ \\ 
 Code error estimate (sup norm): $4.7 \times 10^{-6}$ \\ \hline
 Determinant of residue matrix: $2.3\times 10^{-13}-i5.7\times 10^{-18}$  
 \\ \hline
 Value of $|a_1|$: $-9.7\times 10^{-3}$. \\
 Integration tolerance for $\Psi$: $10^{-9}$ \\ \hline \hline
\end{tabular} 
\caption{Results for Equation 2}
\label{table:2}
\end{table}

\item[Equation 3] The differential equation
\be y^{(iv)} - 2(\mbox{e}^xy')' + (\mbox{e}^{2x}-\mbox{e}^x)y = \lambda y 
\label{eq:na2} 
\ee
on the interval $[0,\infty)$. The differential operator here is the formal 
square of the second order operator $\ell y = -y'' + \exp(x)y$. With Dirichlet
conditions $y(0) = 0 = y'(0)$ the eigenvalues are
not known in closed form. However it is a relatively
straightforward matter to compute approximations using
the code SLEUTH \cite{kn:G+M}. For example, computing
at different tolerances and using different truncations of $[0,\infty)$,
the following approximations seem to be correct to all decimal places
quoted:
\[ \lambda_0 = 35.560604, \;\;
   \lambda_1 = 128.113477, \;\;
   \lambda_2 = 297.84692.
\]
For Neumann boundary conditions $y'''(0) = 0 = y''(0)$ the
corresponding approximate eigenvalues obtained were
\[ \lambda_0 = 6.199245, \;\;
   \lambda_1 = 43.002631, \;\;
   \lambda_2 = 136.295990. 
\]
>From this numerical evidence there is no overlap between the first few
eigenvalues of the Dirichlet and Neumann spectra. To investigate whether 
or not there are likely to be HELP inequalities associated with this equation, 
we must compute the residues of the Titchmarsh-Weyl matrices at these 
eigenvalues using our code. The results are shown in Table \ref{table:3}. 
These residue matrices appear (to within the error we expected at the given 
tolerance) to be of rank 1. This suggests that there is no HELP inequality
associated with (\ref{eq:na2}).
\begin{table}
\hspace{0.15in}\begin{tabular}{|c|} \hline
 Using $\alpha = 1$, truncating $[0,\infty)$ to $[0,10]$: \\ \hline
 $\mbox{Res}(M_D,\lambda=35.560604) = 
   \left(\begin{array}{cc} -297.2883 & -110.8968  \\
                           -110.8968  & -41.3676  \end{array}\right) + 
   i \left(\begin{array}{cc} -1.8\times 10^{-6} & -6.9\times 10^{-6} \\
                           -6.9\times 10^{-6} & -2.6\times 10^{-7} 
 \end{array}\right)$ \\ 
 Code error estimate (sup norm): $2.4 \times 10^{-4}$ \\ \hline
 Determinant of residue matrix:  $-7.2\times 10^{-5}-i2.0\times 10^{-8}$  
 \\ \hline
 Value of $|a_1|$: $1.3\times 10^{-2}$. \\
 Integration tolerance for $\Psi$: $10^{-9}$ \\ \hline 
 $\mbox{Res}(M_N,\lambda=35.560604) = 
   \left(\begin{array}{cc} -297.2883 & -110.8968 \\
                           -110.8968  & -41.3676  \end{array}\right) + 
   i \left(\begin{array}{cc} -5.7\times 10^{-5} & -2.1\times 10^{-5} \\
                           -2.1\times 10^{-5} & -7.9\times 10^{-6} 
 \end{array}\right)$ \\ 
 Code error estimate (sup norm): $4.0 \times 10^{-4}$ \\ \hline
 Determinant of residue matrix: $-7.3\times 10^{-5}-i1.8\times 10^{-8}$  
 \\ \hline
 Value of  $|a_1|$: $1.3\times 10^{-2}$. \\
 Integration tolerance for $\Psi$: $10^{-9}$ \\ \hline \hline
\end{tabular} 
\caption{Results for Equation 3}
\label{table:3} 
\end{table}

\end{description}

\section{Acknowledgments}
 We would like to thank both referees for their exceptionally careful
 reading of our first draft, which led to substantial improvements.

\end{document}